\mag=\magstep1
\documentstyle{amsppt}

\topmatter
\title On the cone of divisors of Calabi-Yau fiber spaces
\endtitle
\author Yujiro Kawamata
\endauthor

\rightheadtext{Cone of divisors of Calabi-Yau fiber spaces}

\address Department of Mathematical Sciences, University of Tokyo, Komaba,
Meguro, Tokyo, 153, Japan \endaddress
\email kawamata\@ms.u-tokyo.ac.jp\endemail



\endtopmatter

\document

\head Introduction 
\endhead

Let $f: X \to S$ be a projective surjective morphism 
of normal varieties with geometrically connected fibers.  
We call it a {\it Calabi-Yau fiber space} if 
$X$ has only $\Bbb Q$-factorial terminal singularities and 
the canonical divisor $K_X$ is relatively numerically trivial over $S$.
This concept is a natural generalization of that of Calabi-Yau manifolds.
Such fiber spaces appear as the output of 
the minimal model program (MMP).
We shall investigate divisors on them by using the log minimal model program
(log MMP).
We refer the reader to [KMM] for the generalities of the minimal model theory.

We shall consider the following generalizations of conjectures of 
D. Morrison ([M1, M2]) concerning the finiteness properties of the cones
which are generated by nef divisors or movable divisors (cf. Definition 1.1):

(1) The number of the $\text{Aut}(X/S)$-equivalence classes of
faces of the {\it effective nef cone} $\Cal A^e(X/S)$ corresponding to
birational contractions or fiber space structures is finite. 
Moreover, there exists a finite rational polyhedral cone $\Pi$ which is a
fundamental domain for the action of $\text{Aut}(X/S)$ on
$\Cal A^e(X/S)$ in the sense that
\roster
\item "(a)" $\Cal A^e(X/S) = 
\bigcup_{\theta \in \text{Aut}(X/S)} \theta_*\Pi$, 

\item "(b)" $\text{Int }\Pi \cap \theta_*\text{Int }\Pi = \emptyset$
unless $\theta_* = \text{id}$.
\endroster

(2) The number of the $\text{Bir}(X/S)$-equivalence classes of
chambers $\Cal A^e(X'/S, \alpha)$ in the {\it effective movable cone}
$\Cal M^e(X/S)$ for the marked minimal models $f': X' \to S$ of
$f:X \to S$ with markings $\alpha: X' -\to X$ is finite.  In other words, 
the number of isomorphism classes of the
minimal models of $f: X \to S$ is finite.
Moreover, there exists a finite rational 
polyhedral cone $\Pi'$ which is a 
fundamental domain for the action of $\text{Bir}(X/S)$ on
$\Cal M^e(X/S)$.

A {\it marked minimal model} is a pair consisting of a minimal model and 
a marking birational map to a fixed model (Definition 1.4).

The {\it nef cone} $\bar \Cal A(X/S)$ is known to be locally rational 
polyhedral inside the {\it big cone} $\Cal B(X/S)$ ([K2], Theorem 1.9).
In the case $\text{dim X} = 3$,
we shall prove a similar statement for 
the {\it movable cone} $\bar \Cal M(X/S)$:
the decomposition of the movable cone into nef cones is locally finite 
inside the big cone (Theorem 2.6), and the {\it pseudo-effective cone}
$\bar \Cal B(X/S)$ itself is locally rational polyhedral away from
$\bar \Cal M(X/S)$ (Theorem 2.9).

It is already known that the above conjectures are true if 
$\text{dim }X = \text{dim }S = 3$ ([KM], Theorem 2.5).
The main result of this paper is the proof of the first parts of the 
conjectures (1) and (2) in the case where 
$0 < \text{dim }S < \text{dim }X = 3$ (Theorems 3.6 and 4.4).
In particular, the number of minimal models in a fixed birational class
of 3-folds is finite up to isomorphisms if the Kodaira dimension is positive
(Theorem 4.5).  

In the course of the proof, 
we shall use the $\Bbb R$-divisors in an essential way.  
In fact, $\Bbb R$-divisors are more suitable for the analysis of 
the infinity than the $\Bbb Q$-divisors.

It is necessary to consider the birational version (2) of the conjecture 
in order to carry out our proof for the biregular version (1) 
(cf. Lemmas 1.15 and 1.16).

The relative setting over the base space $S$ is also essential in our 
inductive argument on $\text{dim }S$ with fixed $\text{dim }X$.
This relative setting seems to correspond to the geometric situation where
the size of the metric of the base space $S$ goes to infinity.

\head 1.  Ample cone and movable cone
\endhead

\definition{Definition 1.1}
In this paper, $f: X \to S$ will always be a projective surjective morphism 
of normal varieties defined over $\Bbb C$
with geometrically connected fibers unless stated otherwise.  
A Cartier divisor $D$ on $X$ is said to be
{\it $f$-nef} (resp. {\it $f$-movable}, {\it $f$-effective}, {\it $f$-big})
if $(D \cdot C) \ge 0$ holds for any curve $C$ on $X$ 
which is mapped to a point on $S$
(resp. if $\text{dim Supp Coker}(f^*f_*\Cal O_X(D) \to 
\Cal O_X(D)) \ge 2$, if $f_*\Cal O_X(D) \ne 0$, 
if $\kappa(X_{\eta}, D_{\eta}) = \text{dim }X - \text{dim }S$
for the generic point $\eta \in S$).
A linear combination of Cartier divisors with coefficients in $\Bbb R$ is
called an {\it $\Bbb R$-Cartier divisor}.
The real vector space
$$
\align
N^1(X/S) &= \{\text{Cartier divisor on }X\}/(\text{numerical equivalence
over }S) \otimes_{\Bbb Z} \Bbb R \\
&= \{\Bbb R\text{-Cartier divisor on }X\}/(\text{numerical equivalence over }S)
\endalign
$$
is finite dimensional.  We set 
$\rho(X/S) = \text{dim }N^1(X/S)$.
The class of an $\Bbb R$-Cartier divisor $D$ in $N^1(X/S)$ is denoted by $[D]$.

The {\it $f$-nef cone} $\bar \Cal A(X/S)$ 
(resp. the {\it closed $f$-movable cone} $\bar \Cal M(X/S)$,
the {\it $f$-pseudo-effective cone} $\bar \Cal B(X/S)$) 
is the closed convex cone in $N^1(X/S)$
generated by the numerical classes of $f$-nef divisors 
(resp. $f$-movable divisors, $f$-effective divisors).
We have the following inclusions:
$$
\bar \Cal A(X/S) \subset
\bar \Cal M(X/S) \subset
\bar \Cal B(X/S) \subset N^1(X/S)
$$
The interior $\Cal A(X/S) \subset \bar \Cal A(X/S)$ 
(resp. $\Cal B(X/S) \subset \bar \Cal B(X/S)$) 
is the open convex cone generated by the numerical classes of 
$f$-ample divisors ($f$-big divisors) and called 
an {\it $f$-ample cone} (resp. {\it $f$-big cone}).
We do not know such a characterization for the interior of 
$\bar \Cal M(X/S)$.
We denote by $\Cal B^e(X/S)$ the {\it $f$-effective cone},
the convex cone generated by $f$-effective Cartier divisors. 
We call $\Cal A^e(X/S) = \bar \Cal A(X/S) \cap \Cal B^e(X/S)$ and
$\Cal M^e(X/S) = \bar \Cal M(X/S) \cap \Cal B^e(X/S)$ the 
{\it $f$-effective $f$-nef cone} and {\it $f$-effective $f$-movable cone}, 
respectively. 
By definition, we have $\Cal B(X/S) \subset \Cal B^e(X/S)$.
\enddefinition

\definition{Remark 1.2}
(1) The base space $S$ can be a complex analytic space if we make a suitable
modification.  
In this case, one needs an additional assumption which guarantees 
the finiteness of $\rho(X/S)$.
For example, we consider $S$ as a germ of a neighborhood of a compact
subset $\Cal K \subset S$ as in [K2].

(2) If the log abundance theorem for $\Bbb R$-divisors holds, e.g., 
if $\text{dim }X = 3$ ([KeMM] and [Sho]), 
then $\Cal A^e(X/S)$ and $\Cal M^e(X/S)$ are generated 
by the classes of 
$\Bbb Q$-divisors as convex cones for a Calabi-Yau fiber space $f: X \to S$ 
(Proposition 2.4).
But there may exist rational points in $\bar \Cal M(X/S)$
which do not belong to $\Cal M^e(X/S)$. 
(cf. Example 3.8 (2)).

(3) Even if $[D]$ and $[D']$ are $f$-effective and non-zero, we might
have $[D + D'] = 0$.
Unlike the case of $\bar \Cal A(X/S)$, 
$\bar \Cal M(X/S)$ may contain a linear subspace of $N^1(X/S)$.
\enddefinition

\definition{Example 1.3} 
(1)  If $f: X \to S$ is a birational morphism, then
$\Cal B(X/S) = N^1(X/S)$.

(2) If a generic fiber $X_{\eta}$ of $f: X \to S$ is a curve, 
then the {\it degree} of a divisor $D$ is defined by 
$\text{deg }D = \text{deg }D_{\eta} = (D \cdot F)$ for a general fiber $F$, 
and
$$
\bar \Cal B(X/S) = \{z \in N^1(X/S); \text{deg }z \ge 0\}.
$$ 
\enddefinition

\definition{Definition 1.4}
A {\it minimal model} of $f: X \to S$ (or of $X$ over $S$) is
a projective morphism $f': X' \to S$ which satisfies the following conditions
(cf. [KMM]):

\roster
\item There exists a birational map $\alpha: X' -\to X$ such that
$f' = f \circ \alpha$. 

\item  $X'$ has only $\Bbb Q$-factorial terminal singularities. 

\item $K_{X'}$ is $f'$-nef.
\endroster
The pair $(X', \alpha)$ is called a {\it marked minimal model}
with a marking $\alpha$.
If $f: X \to S$ is also minimal, then $\alpha$ is 
an isomorphism in codimension 1, and we obtain an isomorphism 
$\alpha_*: N^1(X'/S) \to N^1(X/S)$ such that
$\alpha_*(\bar \Cal M(X'/S)) = \bar \Cal M(X/S)$ and
$\alpha_*(\bar \Cal B(X'/S)) = \bar \Cal B(X/S)$.
We denote $\alpha_*(\bar \Cal A(X'/S)) = \bar \Cal A(X'/S, \alpha)$.
We sometimes write $\bar \Cal A(X'/S)$ instead of $\bar \Cal A(X'/S, \alpha)$
if there is no danger of confusion.
Two marked minimal models $(X_i, \alpha_i)$ $(i = 1,2)$
are said to be {\it isomorphic} if
there exists an isomorphism $\beta: X_1 \to X_2$
such that $\alpha_1 = \alpha_2 \circ \beta$.

Let $f: X \to S$ be a minimal model.
We denote by $\text{Aut}(X/S)$ (resp. $\text{Bir}(X/S)$) the group of biregular
(resp. birational) automorphisms of $X$ over $S$.
Any $\theta \in \text{Bir}(X/S)$ doen not contract any divisor on $X$, 
since $K_X$ is $f$-nef and $X$ has only terminal singularities.   
Thus there is a linear representation
$$
\sigma: \text{Bir}(X/S) \to GL(N^1(X/S), \Bbb Z)
$$
given by $\sigma(\theta)([D]) = \theta_*([D])$. 
\enddefinition

\proclaim{Lemma 1.5}
Let $(X_i, \alpha_i)$ $(i = 1,2)$ be marked minimal models of a minimal
$f: X_0 \to S$.
Then the following conditions are equivalent:

(1) $(X_1, \alpha_1)$ and $(X_2, \alpha_2)$ are isomorphic.

(2) $\Cal A(X_1/S, \alpha_1) = \Cal A(X_2/S, \alpha_2)$ in $N^1(X/S)$.

(3) $\Cal A(X_1/S, \alpha_1) \cap \Cal A(X_2/S, \alpha_2) \ne \emptyset$
in $N^1(X/S)$.
\endproclaim

\demo{Proof}
Assume that there exists a divisor $D$ such that
$\alpha_{1*}^{-1}D$ and $\alpha_{2*}^{-1}D$ are relatively ample over $S$.
Then 
$$
X_1 = \text{Proj}(\bigoplus_{m\ge 0} f_{1*}\Cal O_{X_1}(m\alpha_{1*}^{-1}D))
\cong \text{Proj}(\bigoplus_{m\ge 0} f_{2*}\Cal O_{X_2}(m\alpha_{2*}^{-1}D))
= X_2
$$ 
where the isomorphism is compatible with the $\alpha_i$.
Thus (3) implies (1).  Other implications are obvious.
\hfill $\square$
\enddemo

\proclaim{Corollary 1.6}
There is a 1-1 correspondence between the orbit space
$$
\text{Bir}(X_0/S)/\text{Aut}(X_0/S)
$$ 
and the set of isomorphism classes of the
marked minimal models $(X, \alpha)$ of $X_0$ over $S$ such that
$X$ is isomorphic to $X_0$ over $S$.
\hfill $\square$
\endproclaim

\definition{Definition 1.7}
$f: X \to S$ is said to be a {\it Calabi-Yau fiber space} if 
$X$ has only $\Bbb Q$-factorial terminal singularities and 
$[K_X] = 0$ in $N^1(X/S)$.
This concept is more general than the usual Calabi-Yau manifold in the 
following points: (1) there is no assumption on the fundamental group nor
the irregularity of the generic fiber, 
(2) $X$ may be mildly singular, (3) we consider relatively
over the base space $S$. 
For example, if $\text{dim }X = \text{dim }S$ (resp. 
$= \text{dim }S + 1$), then 
$f$ is a {\it crepant resolution} of singularities
(resp. an elliptic fibration).
We note that $h^1(\Cal O_{X_{\eta}})$ may be non-zero even if $h^1(\Cal O_X)
= 0$.
Any minimal model which satisfies the abundance theorem yields a 
Calabi-Yau fiber space (cf. [KMM]).
The point is that we can treat these cases in a unified way.
\enddefinition

\definition{Definition 1.8}
Let $f: X \to S$ be a Calabi-Yau fiber space, and
$X @>g>> T @>h>> S$ a factorization such that $g$ is also a Calabi-Yau 
fiber space and $h$ is not an isomorphism.
Then $g^*: N^1(T/S) \to N^1(X/S)$ is injective, and 
$g^*\bar \Cal A(T/S) = g^*N^1(T/S) \cap \bar \Cal A(X/S)$
is a face of $\bar \Cal A(X/S)$.
There are 2 cases:

(1) $g$ is a birational morphism.  In this case, it is called a 
{\it birational contraction}.  
We have $\rho(X/T) + \rho(T/S) = \rho(X/S)$ ([KMM]). 
If $\rho(X/T) = 1$, then
it is called {\it elementary} or {\it primitive}.

(2) $\text{dim }X > \text{dim }T$.  In this case, $g$ is called a 
{\it fiber space structure} (cf. [O]).

Let $D$ be an $f$-effective but not $f$-nef $\Bbb R$-divisor.
If $\epsilon$ is a sufficiently small positive number, then the pair 
$(X, \epsilon D)$ is log terminal, and there exists an extremal ray $R$ for 
this pair ([KMM]).  
Let $\phi: X \to Y$ be a contraction morphism over $S$ associated to $R$.
Since $K_X$ is $f$-nef, $\phi$ is a primitive birational contraction morphism.
It is called a {\it divisorial contraction} or a {\it small contraction}
if the exceptional locus of $\phi$ is a prime divisor or not, respectively.
In the latter case, the log flip of $\phi$ is called a {\it $D$-flop}.

In this paper, a prime divisor $E$ on $X$ is said to be {\it $f$-exceptional} 
if there exists a minimal model $f': X' \to S$ of $f$ 
and a divisorial contraction 
$\phi: X' \to Y$ over $S$ whose exceptional divisor is the strict transform 
of $E$.
\enddefinition

As a consequence of the cone theorem ([KMM]), we obtain

\proclaim{Theorem 1.9} ([K2, Theorem 5.7])
Let $f: X \to S$ be a Calabi-Yau fiber space.
Then the cone 
$$
\bar \Cal A(X/S) \cap \Cal B(X/S) = \Cal A^e(X/S) \cap \Cal B(X/S)
$$
is locally rational polyhedral inside the open cone $\Cal B(X/S)$. 
Moreover, any face $F$ of this cone 
corresponds to a birational contraction $\phi: X \to Y$ over $S$
by the equality $F = \phi^*(\bar \Cal A(Y/S) \cap \Cal B(Y/S))$.
\hfill $\square$
\endproclaim

The following is an easy generalization of the characterization of 
nef and big divisors in [K1, Lemma 3] to $\Bbb R$-divisors:

\proclaim{Proposition 1.10}
Let $f: X \to S$ be a proper morphism of normal varieties.  Then
$$
\bar \Cal A(X/S) \cap \Cal B(X/S) = 
\{z \in \bar \Cal A(X/S); z_{\eta}^n > 0\}
$$
where $n$ is the dimension of the generic fiber $X_{\eta}$ of $f$.
\hfill $\square$
\endproclaim

The following is in [W] for the case of Calabi-Yau 3-folds:

\proclaim{Corollary 1.11} 
Let $f: X \to S$ be a Calabi-Yau fiber space, and
let $\Cal W = \{z \in N^1(X/S); z_{\eta}^n > 0\}$.
Then the cone $\bar \Cal A(X/S) \cap \Cal W$
is locally rational polyhedral inside the cone $\Cal W$. 
\hfill $\square$
\endproclaim

Now we consider a generalization of the Morrison conjecture:

\proclaim{Conjecture 1.12} (cf. [M1, M2]).
Let $f: X \to S$ be a Calabi-Yau fiber space.  Then the following hold:

(1) The number of the $\text{Aut}(X/S)$-equivalence classes of
faces of the cone $\Cal A^e(X/S)$ corresponding to
birational contractions or fiber space structures is finite. 
Moreover, there exists a finite rational polyhedral cone $\Pi$ which is a
fundamental domain for the action of $\text{Aut}(X/S)$ on
$\Cal A^e(X/S)$ in the sense that
\roster
\item "(a)" $\Cal A^e(X/S) = \bigcup_{\theta \in \text{Aut}(X/S)} 
\theta_*\Pi$, 

\item "(b)" $\text{Int }\Pi \cap \theta_*\text{Int }\Pi = \emptyset$
unless $\theta_* = \text{id}$.
\endroster

(2) The number of the $\text{Bir}(X/S)$-equivalence classes of
chambers $\Cal A^e(X'/S, \alpha)$ in the cone $\Cal M^e(X/S)$ 
for the marked minimal models $f': X' \to S$ of
$f$ is finite.  In other words, 
the number of isomorphism classes of the
minimal models of $f$ is finite.
Moreover, there exists a finite rational 
polyhedral cone $\Pi'$ which is a 
fundamental domain for the action of $\text{Bir}(X/S)$ on
$\Cal M^e(X/S)$.
\hfill $\square$
\endproclaim

\definition{Remark 1.13}
(1) The above conjectures were inspired by the mirror symmetry conjecture of 
Calabi-Yau threefolds.  Some positive evidences are given in [B], [GM] and [OP]
for (1) and [Nm1] for (2).

(2) With respect to our relative formulation over the base space $S$, 
the variety $X$ can be an arbitrary minimal model which satisfies the 
abundance theorem ([KMM]), if we take $S$ to be
the canonical model $\text{Proj}(\bigoplus_{m=0}^{\infty}H^0(X, mK_X))$.

(3) If we replace the ample cone $\Cal A(X)$ by the K\"ahler cone 
$\Cal K(X)$, then the conjecture is clearly false.

(4) The finiteness questions such as 
the finite generation of the canonical ring, the termination of flips,
the finiteness of the cones, the boundedness of the moduli space
and the Zariski decomposition, seem to be mutually related (cf. [A], [G]).
\enddefinition

\definition{Example 1.14}
(1) Let $X$ be an abelian variety.  
Then we have 
$$
\bar \Cal A(X) = \bar \Cal B(X) = 
\{z \in N^1(X) ; z^n \ge 0\}^o
$$
where $^o$ denotes an irreducible component of the cone,
since $X$ does not contain a rational curve,
and there is no divisorial contraction nor flop of $X$.

Although the shape of this cone is quite different from a finite rational 
polyhedral cone, the conjecture seems to be true in this case, too.
One checks it by an explicit calculation 
in the case where $X \cong E \times \cdots \times E$ for an elliptic curve $E$
without complex multiplications (Corollary 2.11).
A related result is in [NN].

(2) Let $X$ be a K3 surface with an ample class $h$, 
and $\Sigma$ the set of all the $(-2)$-curves
on $X$. Then
$$
\bar \Cal A(X) = \{z \in N^1(X) ; z^2 \ge 0, z \cdot h \ge 0, 
z \cdot C \ge 0 \quad \forall C \in \Sigma\} 
$$
and $\bar \Cal B(X)$ is the closed convex cone generated by the cone
$\{z \in N^1(X) ; z^2 \ge 0, z \cdot h \ge 0\}$ and the $C \in \Sigma$.
This duality between $\bar \Cal A$ and $\bar \Cal B$
will be generalized in Theorem 2.9.  
In this case, the conjecture is verified in [St] (see Theorem 2.1). 
See also [Kov].
\enddefinition

Our strategy is to analyse the birational automorphism group first and try
to prove Conjecture (2), and then consider the biregular automorphism
group toward Conjecture (1).

\proclaim{Lemma 1.15}
Let $f: X \to S$ be a Calabi-Yau fiber space.
Assume that the number of faces of $\Cal A^e(X/S)$ 
which correspond to primitive birational contractions is finite up
to the action of $\text{Bir}(X/S)$.
Then it is also finite up to the action of $\text{Aut}(X/S)$.
\endproclaim

\demo{Proof}
Let $F$ be a face of $\Cal A^e(X/S)$ which corresponds to a 
primitive birational contraction.
Let $\theta$ and $\theta'$ be birational automorphisms of $f$ such that
$\theta_*F$ and $\theta'_*F$ are also faces of $\Cal A^e(X/S)$.
Assume that $\theta, \theta' \not\in \text{Aut}(X/S)$.
Then $\theta$ and $\theta'$ map the interior of $\Cal A^e(X/S)$ near $F$ 
to the exterior of it near $\theta_*F$ and $\theta'_*F$, respectively.
Then $\theta' \circ \theta^{-1}$ maps the interior of $\Cal A^e(X/S)$ near 
$\theta_*F$ to the interior of it near $\theta'_*F$, hence
$\theta' \circ \theta^{-1} \in \text{Aut}(X/S)$, because
$$
\text{Aut}(X/S) = \{\theta \in \text{Bir}(X/S) ;
\theta_*\Cal A(X/S) \cap \Cal A(X/S) \ne \emptyset\}
$$
Thus $\theta_*F$ and $\theta'_*F$ belong to the
same $\text{Aut}(X/S)$-orbit of the faces.
Therefore, a $\text{Bir}(X/S)$-orbit splits into at most two
$\text{Aut}(X/S)$-orbits.
\hfill $\square$
\enddemo

\proclaim{Lemma 1.16}
Let $f: X \to S$ be a Calabi-Yau fiber space.
Assume that the number of faces of $\Cal A^e(X/S)$ 
which correspond to fiber space structues is finite up
to the action of $\text{Bir}(X/S)$.
In addition, assume that the first part of 
Conjecture 1.12 (2) is true for any Calabi-Yau
fiber space which factors $f$ non-trivially.
Then it is also finite up to the action of $\text{Aut}(X/S)$.
\endproclaim

\demo{Proof}
Let $F$ be a face of $\Cal A^e(X/S)$ which corresponds to a 
fiber space structure $g: X \to T$.
Let $\theta_{\lambda}$ be birational automorphisms of $f$ such that
the $\theta_{\lambda*}F$ are also faces of $\Cal A^e(X/S)$.
Then $\theta_{\lambda*}^{-1}$ maps the interior of $\Cal A^e(X/S)$ near 
$\theta_{\lambda*}F$ to the interior of some chamber 
$\Cal A^e(X_{\lambda}/S, \alpha_{\lambda})$ near its face $F$.
By the assumption, there exists only finitely many chambers which
have $F$ as faces up to $\text{Bir}(X/T) \cong \text{Bir}(X/S)$.
If $\theta_{\lambda_1*}^{-1}\Cal A^e(X/S) 
= \theta'_*\circ \theta_{\lambda_2*}^{-1}\Cal A^e(X/S)$ for some 
$\theta' \in \text{Bir}(X/T)$, then 
$\zeta = \theta_{\lambda_1}\circ \theta' \circ \theta_{\lambda_2}^{-1} 
\in \text{Aut}(X/S)$ and $\zeta_*\theta_{\lambda_2*}F = \theta_{\lambda_1*}F$,
because $\theta'_*F = F$.
Therefore, a $\text{Bir}(X/S)$-orbit splits into finite
$\text{Aut}(X/S)$-orbits.
\hfill $\square$
\enddemo

\head 2.  General results for dimension 2 or 3.
\endhead

\proclaim{Theorem 2.1}
Let $f: X \to S$ be a Calabi-Yau fiber space such that $\text{dim }X = 2$.
Then Conjecture 1.12 is true.
\endproclaim

\demo{Proof}
Since the minimal models of dimension 2 are absolutely minimal, we have
$\text{Aut}(X/S) = \text{Bir}(X/S)$.
If $S$ is not a point, then the assertion is trivial.
Therefore, we consider the case where $X$ is a K3 surface, an Enriques surface,
an abelian surface, or a hyperelliptic surface over a point $S$.
The first case is proved in [St], and other cases are similar.
We recall briefly the proof for the convenience of the reader.

Let $O(N^1(X))$ be the orthogonal group of $N^1(X)$ with respect to the 
intersection pairing, $\Cal C = \{z \in N^1(X) ; z^2 > 0\}^o$ 
the positive cone, and $N = N^1(X)_{\Bbb Z}$.
We consider the arithmetic subgroup 
$$
\Gamma(X) = \{\theta \in O(N^1(X))_{\Bbb Z} ; \theta(\Cal C) = \Cal C,
\theta \text{ induces the identity on }N^*/N\}.
$$
Then we can find a rational polyhedral fundamental domain $\Pi \subset 
\Cal A^e(X)$
for the action of $\Gamma(X)$ on the convex hull of 
$\bar \Cal C \cap N^1(X)_{\Bbb Q}$ by using a result of [AMRT, Chap. II].
Let $\Gamma_B(X)$ be the subgroup of $\Gamma(X)$ 
which consists of elements preserving nodal classes. 
Then $\Pi$ serves also as a fundamental domain for the action of 
$\Gamma_B(X)$ on $\Cal A^e(X)$.

First, assume that $X$ is a K3 surface.
By the Torelli theorem, elements of $\Gamma_B(X)$ 
correspond to automorphisms of
$X$ which preserve a non-zero holomorphic 2-form on $X$.
Therefore, $\Gamma_B(X)$ can be identified as a subgroup of 
$\text{Aut}(X)$ of finite index, and we have the assertion.

Next, assume that $X$ is an Enriques surface.
Then the universal cover $\tilde X$ is a K3 surface with an involution 
$\sigma$ such that $X = \tilde X/\sigma$.
We have $\text{Aut}(X) \cong \{\theta \in \text{Aut}(\tilde X) ;
\theta \circ \sigma = \sigma \circ \theta\}/\{1, \sigma \}$.
We can extend any automorphism of $N^1(X)$ to that of 
$N^1(\tilde X)$ so that it coincides with 
the identity on the orthogonal complement of $N^1(X)$ in $N^1(\tilde X)$
with respect to the intersection pairing.
Thus elements of $\Gamma_B(X)$ can be extended to automorphisms of $X$, and we 
have our assertion.

If $X$ is an abelian or a hyperelliptic surface, 
then we replace $\Gamma(X)$ by its subgroup
consisting of elements whose determinant equal 1 ([Shi]).  
Then the rest is the same.
\hfill $\square$
\enddemo

\definition{Remark 2.2}
The above theorem is also valid over any field $k$ of characteristic 0.
Indeed, if $X$ is defined over $k$ and $\bar X = X \times_k \bar k$,
then $N^1(X/k) = N^1(\bar X/\bar k)^G$ 
for the Galois group $G = \text{Gal}(\bar k/k)$.
Assume that $X$ is a K3 surface (other cases are similar).
We can extend any element of the similarly defined $\Gamma(X)$ to
$\Gamma(\bar X)$ such that it acts trivially on the orthogonal complement
of $N^1(X/k)$ in $N^1(\bar X/\bar k)$
with respect to the intersection pairing.
If $\theta$ is the corresponding automorphism of $\bar X$, then we have
$\sigma \circ \theta \circ \sigma^{-1}= \theta$ in $\text{Aut}(\bar X/\bar k)$
for any $\sigma \in G$, because they induce 
the same automorphism of $H^2(\bar X, \Bbb Z)$.
Hence $\theta \in \text{Aut}(X/k)$, and we obtain a rational polyhedral 
fundamental domain for the action of $\text{Aut}(X/k)$ on 
$\Cal A^e(X/k)$.
\enddefinition

\proclaim{Theorem 2.3} (cf. [K2, p.120]).
Let $f_0: X_0 \to S$ be a Calabi-Yau fiber space with $\text{dim }X_0 = 3$, and
$D$ an $\Bbb R$-divisor such that $[D] \in \Cal M^e(X_0/S)$.
Then there exists a sequence of $D$-flops such that the strict transform of 
$D$ becomes relatively nef over $S$.
Therefore, 
$$
\Cal M^e(X_0/S) = \bigcup_{(X, \alpha)} \Cal A^e(X/S, \alpha)  
$$
where the union on the right hand side is taken for all the marked minimal 
models $(X, \alpha)$ of $X_0$ over $S$. 
\endproclaim 

\demo{Proof}
If $D$ is not $f_0$-nef, then there exists an $f_0$-ample $\Bbb R$-divisor 
$D_1$ such that $D + D_1$ becomes a $\Bbb Q$-divisor which is not $f_0$-nef but
$[D+D_1] \in \Cal M^e(X_0/S)$.
Then there exists a $D$-flop by [K2, 2.4].
The proof of the termination of any sequence of $D$-flops in [K3, Lemma 4] 
works also for $\Bbb R$-divisors.
\hfill $\square$
\enddemo

\proclaim{Proposition 2.4}
Let $f: X \to S$ be a Calabi-Yau fiber space such that $\text{dim }X = 3$. 
Then the cones $\Cal A^e(X/S)$ and $\Cal M^e(X/S)$ are generated by 
the numerical classes of $\Bbb Q$-Cartier divisors as convex cones.
\endproclaim 

\demo{Proof}
By the log abundance theorem, any $f$-effective $f$-nef
$\Bbb R$-divisor is a pull-back of an ample
$\Bbb R$-Cartier divisor by some morphism ([Sho]).
Since the ample cone is generated by the classes of Cartier divisors, we
have the assertion for $\Cal A^e(X/S)$.
Any $f$-effective $f$-movable $\Bbb R$-divisor becomes
relatively nef after a finite number of flops.
Hence the assertion for $\Cal M^e(X/S)$.
\hfill $\square$
\enddemo

The following gives a positive answer to Conjecture 1.12
in a special case, where we note that $\text{Bir}(X/S) = \{\text{id}\}$:

\proclaim{Theorem 2.5} ([KM]). 
Let $S$ be a normal 3-fold, and $f: X \to S$ a minimal resolution.
Then $\bar \Cal A(X/S)$ is a finite polyhedral cone, and there
exists only finitely many marked minimal models of $f$.
In other words, Conjecture 1.12 is true if $\text{dim }X = \text{dim }S = 3$.
\hfill $\square$
\endproclaim

The following is a generalization of the above theorem:

\proclaim{Theorem 2.6}
Let $f_0: X_0 \to S$ be a Calabi-Yau fiber space 
such that $\text{dim }X_0 = 3$. 
Then the decomposition 
$$
\Cal M^e(X_0/S) \cap \Cal B(X_0/S) = 
\bigcup_{(X, \alpha)} \Cal A^e(X/S, \alpha) \cap \Cal B(X_0/S)
$$
is locally finite inside the open cone $\Cal B(X_0/S)$ in the following sense:
if $\Sigma$ is a closed convex cone contained in 
$\Cal B(X_0/S) \cup \{0\}$, then there exist only a finite number of
cones $\Cal A^e(X/S, \alpha) \cap \Cal B(X_0/S)$ which intersect $\Sigma$.
\endproclaim 

\demo{Proof}
Let $D$ be any $\Bbb R$-divisor such that 
$[D] \in \Cal M^e(X_0/S) \cap \Cal B(X_0/S)$.
Then there exists a marked minimal model $(X_1, \alpha_1)$ such that
$[D] \in \Cal A^e(X_1/S, \alpha_1)$.
Let $F$ be the face of $\Cal A^e(X_1/S, \alpha_1)$ whose interior
contains $[D]$.  Then there exists
a projective birational contraction morphism $\phi: X_1 \to Y$ over $S$ to a
variety $Y$ with only canonical singularities such that 
$F = \phi^*\Cal A^e(Y/S)$. 
By [KM], there exists only a finite number of cones
$\Cal A^e(X_i/Y, \beta_i)$ for the marked minimal models 
$(X_i, \beta_i)$ of $X_1$ over $Y$.
This implies that there exist only a finite number of cones
$\Cal A^e(X_i/S, \alpha_1 \circ \beta_i)$ 
for the marked minimal models of $X_0$ over $S$ which
contain $F$ as their faces.
Since these cones are locally rational 
polyhedral in the open cone
$\Cal B(X_0/S)$, there exists a small open cone $\Sigma'$ which contains 
$[D]$ such that $\Cal M^e(X_0/S) \cap \Sigma' \subset 
\bigcup_i \Cal A^e(X_i/S, \alpha_1 \circ \beta_i)$.  
This proves our theorem.
\hfill $\square$
\enddemo

The accumulation occurrs only toward the boundary 
$\partial \bar \Cal B(X/S)$:

\proclaim{Corollary 2.7}
Let $f: X \to S$ be a Calabi-Yau fiber space such that
$\text{dim }X = 3$. 
Then the cone $\bar \Cal M(X/S) \cap \Cal B(X/S)$ is locally rational 
polyhedral inside the open cone $\Cal B(X/S)$.
Moreover, the faces of this cone correspond to divisorial contractions of some
marked minimal models.
\hfill $\square$
\endproclaim

\definition{Remark 2.8}
One $f$-exceptional divisor may correspond to several faces of 
\linebreak
$\bar \Cal M(X/S)$.
\enddefinition

\proclaim{Theorem 2.9}
Let $f: X \to S$ be a Calabi-Yau fiber space such that
$\text{dim }X = 3$. 
Then the cone $\bar \Cal B(X/S)$ is locally rational 
polyhedral inside the open cone 
$$
N^1(X/S) \setminus (\bar \Cal M(X/S) \cap \partial \bar \Cal B(X/S)).
$$
Moreover, it is generated by $\bar \Cal M(X/S)$ and the 
numerical classes of the $f$-exceptional divisors.
\endproclaim

\demo{Proof}
Let $U$ be an open cone which contains
$\bar \Cal M(X/S) \cap \partial \bar \Cal B(X/S) \setminus \{0\}$.
The cone $\bar \Cal M(X/S) \setminus U$ 
inside $N^1(X/S) \setminus U$ is bounded by 
a finite number of faces $F_i$ $(i = 1, \ldots, \ell)$ 
of some of the chambers  
$\bar \Cal A(X_i/S, \alpha_i)$ 
for marked minimal models $(X_i, \alpha_i)$ of $f$, where 
some of the $X_i$ may coincide.
Let $h_i$ $(i = 1, \ldots, \ell)$ be a linear functional
on $N^1(X/S)$ which defines a hypersurface containing $F_i$ and such that
$$
\bar \Cal M(X/S) \setminus U = 
\{z \in N^1(X/S) \setminus U ; h_i(z) \ge 0 \quad \forall i \}
$$
For any $i = i_0$,
there exists only one prime divisor $D = D_{i_0}$ on $X$ 
such that $h_{i_0}(D) < 0$.
Indeed, we have $(\alpha_{i_0*}^{-1}D \cdot C) < 0$ 
for any extremal rational curve $C$ 
of the divisorial contraction of $X_{i_0}$ corresponding to
the face $F_{i_0}$, and $\alpha_{i_0*}^{-1}D$
coincides with the exceptional divisor of this contraction.
Therefore, we have our assertion.
\hfill $\square$
\enddemo

We have the following positive evidence for Conjecture 1.12
in the case where $X$ is a direct product of an elliptic curve without complex
multiplications.  By the mirror symmetry, its complexified K\"ahler cone 
$\Bbb R^{\frac 12 n(n+1)} \times \sqrt{-1}\Cal A(X)$
should be isomorphic to the moduli space
of marked principally polarized abelian varieties under the mirror map,
as is proved in the following proposition.
An abelian variety which is isogenous to $X$ may correspond to 
a non-principal polarization.

\proclaim{Proposition 2.10}
Let $X = E \times \cdots \times E$ ($n$-times) 
for an elliptic curve $E$ without complex multiplications.
Then $\text{Aut}(X) = \text{Bir}(X)$, $\rho(X) = \frac 12 n(n+1)$, and
$$
\text{Im}(\sigma: \text{Aut}(X) \to GL(N^1(X), \Bbb Z)) \cong
GL(n, \Bbb Z).
$$  
Moreover, there is a linear isomorphism $\tau: N^1(X) \to S(n, \Bbb R)$ 
to the real vector space of symmetric $(n,n)$-matrices which sends
$\Cal A(X)$ to the cone of positive definite matrices and
which is compatible with the natural $GL(n, \Bbb Z)$-actions.
\endproclaim

\demo{Proof}
The first part is clear.
We take a basis of $N^1(X)$ which consists of
$D_i = \text{pr}_i^*(0)$ for $1 \le i \le n$
and $E_{ij} = \text{pr}_{ij}^*(\Delta)$ for $1 \le i < j \le n$,
where the $\text{pr}$ are projections to suitable factors and $\Delta$ 
denotes the diagonal.
We claim that
$$
[\sum_{1 \le i \le n} x_{ii}D_i + 
\sum_{1 \le i < j \le n}x_{ij}(D_i + D_j - E_{ij})]^n
= n! \text{ det }\vert x_{ij} \vert,
$$
where we set $x_{ij} = x_{ji}$ if $i > j$.
Indeed, this follows from the following equalities with suitable permutations
of indices:
$$
\align
&D_1D_2\cdots D_n = 1 \\
&D_1^2 = 0 \\
&D_1(D_1+D_2- E_{12}) = 0 \\
&(D_1+D_2- E_{12})^2 = -2D_1D_2 \\
&(D_1+D_2- E_{12})^3 = 0 \\
&(D_1+D_2- E_{12})^2(D_1+D_3- E_{13}) = 0 \\
&(D_1+D_2- E_{12})(D_2+D_3- E_{23})(D_1+D_3- E_{13}) = 2D_1D_2D_3 \\
&(D_1+D_2- E_{12})(D_2+D_3- E_{23})
\cdots (D_{i-1}+D_i- E_{i-1,i})(D_1+D_i- E_{1i}) \\ 
&= (-1)^{i-1}2D_1D_2\cdots D_i 
\endalign
$$
We define $\tau(\sum_{1 \le i \le n} x_{ii}D_i + 
\sum_{1 \le i < j \le n}x_{ij}(D_i + D_j - E_{ij}))
= [x_{ij}]$.
For example, an ample divisor $\sum_i D_i$ on $X$ is mapped to the identity
matrix.  By the above equality,
$D$ is ample if and only if $\tau(D)$ is positive definite for any 
$D \in N^1(X)$.

Let $\theta = (\theta_{ij}) \in GL(n, \Bbb Z)$.  We shall check the equality
$\tau(\theta_*D) = \theta^{*-1}\tau(D)\theta^{-1}$ 
for the following generators of $GL(n, \Bbb Z)$, where $^*$ denotes the
transpose of a matrix:

(1) $\theta_{11} = \theta_{12} = \theta_{22} = \theta_{33} = \cdots = 
\theta_{nn} = 1$ and $\theta_{ij} = 0$ for other $i$ and $j$.

(2) $\theta_{11} = -1$, $\theta_{22} = \theta_{33} = \cdots = \theta_{nn} = 1$ 
and $\theta_{ij} = 0$ for other $i$ and $j$.

(3) $\theta_{12} = \theta_{21} = \theta_{33} = \cdots = \theta_{nn} = 1$ and
$\theta_{ij} = 0$ for other $i$ and $j$.

Case (1).  $\theta_*D_1=E_{12}$, $\theta_*D_i = D_i$ for $i \ge 2$,
$\theta_*E_{12} = -D_1 + 2D_2 + 2E_{12}$,
$\theta_*E_{1i} = -D_1 + D_2 + D_i + E_{12} + E_{1i} - E_{2i}$ for $i \ge 3$,
and $\theta_*E_{ij} = E_{ij}$ for other $i$ and $j$.
If we set $\tau(D) = [x_{ij}]$ and $\tau(\theta_*D) = [y_{ij}]$, 
then $y_{11}= x_{11}$, $y_{22} = x_{11}+x_{22}-2x_{12}$,
$y_{12}=-x_{11} + x_{12}$, $y_{1i}=x_{1i}$, and $y_{2i}=-x_{1i}+x_{2i}$. 
Thus $[y_{ij}] = \theta^{*-1}[x_{ij}]\theta^{-1}$.

Case (2).  $\theta_*D_i = D_i$ for all $i$,
$\theta_*E_{1i} = 2D_1 + 2D_i - E_{1i}$ for $i \ge 2$,
and $\theta_*E_{ij} = E_{ij}$ for other $i$ and $j$.
Then $y_{1i}= -x_{1i}$ for $i \ge 2$ and 
$y_{ij} = x_{ij}$ for other $i$ and $j$.
Thus $[y_{ij}] = \theta^{*-1}[x_{ij}]\theta^{-1}$.                          
Case (3) is clear.
\hfill $\square$
\enddemo

\proclaim{Corollary 2.11}
Conjecture 1.12 is true for $X = E \times \cdots \times E$
for an elliptic curve $E$ without complex multiplications.
\hfill $\square$
\endproclaim

We make a remark on the behaviour of the cones of divisors under deformations
extending [W] and [Nm2].  
This result is not used in the rest of this paper.

\proclaim{Proposition 2.12} 
Let $X$ be a Calabi-Yau fiber space over a point such that
$\text{dim }X = 3$ an $h^2(\Cal O_X) = 0$. 
Let $\pi: \Cal X \to B$ be a flat family of deformations of 
$X = X_0$ over a germ $(B, 0)$.
Then there exist at most countably many proper closed analytic subsets 
$C_{\lambda}$ of $B$, which may contain $0$,
such that $\bar \Cal A(X_t)$, $\bar \Cal M(X_t)$ and 
$\bar \Cal B(X_t)$ are constant in $N^1(X) \cong N^1(X_t)$ for
$t \in B \setminus \bigcup_{\lambda} C_{\lambda}$.
\endproclaim 

\demo{Proof}
Since $h^2(\Cal O_X) = 0$, any line bundle on $X_0$ can be extended to 
$\Cal X$.  Hence $\pi$ is projective,
and $\rho(X_t)$ is constant on $t \in B$ (cf. [Nm2, 6.1]).
By [KoM, 12.1.10], $X_t$ is $\Bbb Q$-factorial for any $t \in B$.  
We claim that any flop $X_0 @>{\phi_0}>> Y_0 @<{\phi_0^+}<< X_0^+$ 
can be extended to a simultaneous flop 
$\Cal X @>{\phi}>> \Cal Y @<{\phi^+}<< \Cal X^+$ over $B$.
Indeed, a free linear system $\vert L_0\vert$ on $X_0$ which gives $\phi_0$
is extended to a free linear system $\vert \Cal L\vert$ to give $\phi$.
The induced morphism $\phi_t$ on $X_t$ is not an isomorphism, because
the exceptional curves of $\phi_0$ deform to $\Cal X$.
The simultaneous flop $\phi^+$ of $\phi$ exists as usual.

Let $\Cal D$ be a prime divisor of $\Cal X$ which does not contain a fiber
of $\pi$ such that $D_0$ is an exceptional prime divisor of $X_0$.
We claim that $D_t$ is exceptional in $X_t$ for all $t \in B$.
In order to prove this, let $X_0 = X_0^0 @>{\phi_0^0}>> 
Y_0^0 @<{\phi_0^{0+}}<< X_0^1 @>{\phi_0^1}>> \cdots @<{\phi_0^{N-1+}}<< X_0^N$ 
be a sequence of $D_0$-flops with a 
divisorial contraction $X_0^N @>{\phi_0^N}>> Y_0^N$ 
whose exceptional divisor is the strict transform $D_0^N$ of $D_0$.
As an extension of the above sequence over $B$, 
we obtain a sequence of simultaneous $\Cal D$-flops 
$\Cal X = \Cal X^0 @>{\phi^0}>> \Cal Y^0 @<{\phi^{0+}}<< \Cal X^1 @>{\phi^1}>> 
\cdots @<{\phi^{N-1+}}<< \Cal X^N$ 
with a primitive birational contraction $\Cal X^N @>{\phi^N}>> \Cal Y^N$,
which should be a divisorial contraction of $\Cal D$ by the dimension count.
Because prime divisors on the $X_t$ form countable number of families over $B$,
we obtain our assertion.
\hfill $\square$
\enddemo

\head 3.  Elliptic fiber spaces
\endhead

We shall consider a Calabi-Yau fiber space $f: X \to S$ such that
$\text{dim }X = 3$ and $\text{dim }S = 2$ in this section,
and prove that the first parts of Conjecture 1.12 (1) and (2) are true
in this case.
A prime divisor $D$ on $X$ is called {\it $f$-vertical} if $f(D) \ne S$. 
An $f$-exceptional divisor is always $f$-vertical.

\proclaim{Lemma 3.1}
Let $D$ be an $f$-vertical prime divisor.  

(1) If $f(D)$ is a point, then $D$ is $f$-exceptional.

(2) If $f(D)$ is a curve and $D$ does not contain the fiber of 
the morphism $f^{-1}f(D) 
\mathbreak
\to f(D)$ over the generic point of $f(D)$, 
then $D$ is $f$-exceptional.

(3) If $f(D)$ is a curve and $D$ contains the fiber of 
the morphism $f^{-1}f(D) \to f(D)$ over the generic point of $f(D)$, 
then $mD$ is $f$-movable for some positive integer $m$.
\endproclaim

\demo{Proof}
Since $D$ never becomes $f$-nef after flops over $S$ in the case (1) or (2),
$D$ is $f$-exceptional.
(3) Since $S$ has only rational singularities by [Kol] or [Nk1],
it is $\Bbb Q$-factorial.
Let $C$ be a Cartier divisor on $S$ whose support is $f(D)$.
Then $f^*C = mD + D'$ for a positive integer $m$ and
an effective divisor $D'$ such that 
$\text{dim }f(D') = 0$.
Then $mD$ is $f$-movable, since 
there exists an effetive Cartier divisor $C'$ on $S$ which does not contain
$f(D)$ and such that $f^*C' \ge D'$.  
\hfill $\square$
\enddemo

\proclaim{Lemma 3.2}
Let $V(X/S)$ be the subspace of $N^1(X/S)$ generated 
by the classes of $f$-vertical divisors.  Then

(1) $V(X/S) = \Cal B^e(X/S) \cap \partial \bar \Cal B(X/S)$.

(2) $V(X/S)$ is generated by the classes of $f$-exceptional divisors as 
an $\Bbb R$-vector space.
\endproclaim

\demo{Proof}
(1) Let $D$ be an $f$-vertical prime divisor.  
Then there exists a prime divisor
$C$ on $S$ which contains $f(D)$.
We have $f^*C = mD + D'$ for a positive rational number 
$m$ and an effective $\Bbb Q$-divisor
$D'$.  Therefore, $-[D] \in \Cal B^e(X/S)$.
Conversely, an $f$-effective $\Bbb R$-divisor of degree 0 is $f$-vertical.
(2) If $D$ is not an $f$-exceptional 
divisor, then $D'$ is a sum of $f$-exceptional divisors.
\hfill $\square$
\enddemo

\proclaim{Lemma 3.3}
(1) There exists a uniquely determined factorization 
$X @>g>> T @>h>> S$ by a Calabi-Yau fiber space $g$ and 
a projective birational morphism $h$
from a normal surface $T$ which is maximal in the sense that
any other such factorization $X @>g_1>> T_1 @>h_1>> S$ is induced from a
further factorization $T \to T_1 \to S$.

(2) There exists a minimal model
$f': X' \to S$ of $f$ with a factorization 
$X' @>{g'}>> T' @>{h'}>> S$ as in (1) such that 
$g'$ has only 1-dimensional fibers.

(3) The marking birational map $\alpha: X' -\to X$ induces a morphism
$\beta: T' \to T$ such that $\beta \circ g' = g \circ \alpha$.

(4) $\bar \Cal A(X'/T')$ is a finite rational polyhedral cone, and
any $g'$-vertical divisor $D$ is $g'$-exceptional if $[D] \ne 0$ in
$N^1(X'/T')$.
\endproclaim

\demo{Proof}
(1) If $D_1$ and $D_2$ are $f$-nef and $f$-effective but not
$f$-big divisors which correspond to some factorizations.
Then $D_1 + D_2$ gives another factorization which dominates both of them. 

(2) We proceed by induction on $\rho(X/S)$.
We shall denote by $f: X \to S$ any minimal model obtained
by flops by abuse of notation.
Suppose that there exists a prime divisor $E$ 
on $X$ such that $f(E)$ is a point.
Then we take a general member $C$ from a non-complete linear system of divisors
on $S$ containing $f(E)$ without a fixed component.
We write $f^*C = D_1 + D_2$, where $D_1$ (resp. $D_2$)
is the sum of irreducible components 
which are mapped to curves (resp. points) on $S$.
Since $D_1$ is movable, there exists a sequence of flops
which makes $D_1$ relatively nef over $S$.
Let $g_1: X \to S_1$ be the morphism corresponding to $D_1$.
Since $D_1$ is not $f$-big, we have $\text{dim }S_1 = 2$.
Since $[D_1] \ne 0$ in $N^1(X/S)$, we have $S_1 \ne S$.
Therefore, $f$ is factorized non-trivially, and we obtain our assertion 
by induction hypothesis.

(3) Let $H$ be an $h$-ample divisor on $T$.
Since $M = \alpha_*^{-1}g^*H$ is $g'$-movable, 
there exists a sequence of flops over $T'$ such that 
the strict transform $M'$ becomes relatively nef over $T'$.
Since $T'$ is maximal, $M'$ is the pull-back of some divisor $H'$ on $T'$.
Then $H'$ is $h'$-movable, hence $h'$-nef, 
and induces a morphism $\beta: T' \to T$.

(4) All the irreducible curves in the fibers of $g'$ form 
finitely many algebraic families.
\hfill $\square$
\enddemo

\proclaim{Corollary 3.4}
There exist only finitely many faces of the chambers in $N^1(X/S)$ which 
correspond to fiber space structures of some marked minimal models of $f$.
\hfill $\square$
\endproclaim

\proclaim{Lemma 3.5}
The image of the representation
$\sigma: \text{Bir}(X/S) \to GL(N^1(X/S), \Bbb Z)$
contains a subgroup $G(X/S)$ of finite index which is isomorphic to a
free abelian group of rank
$\rho(X/S) - v(X/S) - 1$ for $v(X/S) = \text{dim }V(X/S)$.
Moreover, $G(X/S)$ acts properly discontinuously
on an affine linear subquotient space
$W(X/S) = \{z \in N^1(X/S)/V(X/S) ; \text{deg }z = 1\}$ 
as a group of translations such that the quotient space
$W(X/S)/G(X/S)$ is a real torus.
\endproclaim

\demo{Proof}
The degree of divisors is invariant and the space $V(X/S)$ is stable
under the action of $\text{Bir}(X/S)$.
Therefore, we can consider its action on $W(X/S)$.
Moreover, since $V(X/S)$ is generated by the classes of 
$f$-exceptional divisors whose number is finite,
there exists a subgroup of finite index of $\text{Bir}(X/S)$ which
acts trivially on $V(X/S)$.

First, assume that there exists a birational section $D_0$ of $f$.
It corresponds to a $\Bbb C(S)$-rational point of the generic fiber
$X_{\eta}$.
Let $M$ be the group of birational sections of $f$, which is known to be
a finitely generated abelian group.
Since the group of automorphisms of $X_{\eta}$ which fix 
the origin $D_{0, \eta}$
is finite, $M$ is can be identified as a subgroup of finite index
in $\text{Bir}(X/S)$ which acts on $X$ as a group of translations.

Let $\theta$ be a birational automorphism of $f$ which corresponds to
a birational section $D_1$. 
Then $\theta_*D_0 = D_1$.  For a divisor $D$ with $\text{deg }D = 1$, 
we have $\theta_*D_{\eta} - D_{\eta} \sim D_{1,\eta} - D_{0,\eta}$ on
$X_{\eta}$.  Hence $\theta$ acts on $W$ as a translation by
$[D_1 - D_0] \in N^1(X/S)/V(X/S)$. 
The map $\sigma': M \to N^1(X/S)/V(X/S)$ given by 
$\sigma'(D_1) = [D_1 - D_0]$ is a group homomorphism, and its image
is a free abelian group.
For any divisor $D$ on $X$, if we put $D' = D - (\text{deg }D - 1)D_0$,
then $\text{deg }D' = 1$, and 
$f_*\Cal O_X(D')$ is a torsion free sheaf of rank 1 on $S$, hence
$[D']$ is equal to a class of a birational section in $N^1(X/S)/V(X/S)$.
Therefore, the image of $\sigma'$ has rank $\rho(X/S) - v(X/S) -1$.

Now, we consider the general case.
Let $d$ be the smallest positive integer among the degrees of 
the prime divisors on $X$, and fix a prime divisor $D_0$ 
whose degree is $d$.
There exists a subgroup $M$ of finite index of $\text{Bir}(X/S)$ which
acts on $X_{\eta}$ as a group of translations.
For $\theta \in M$ and a divisor $D$ with $\text{deg }D = d$, 
we have again 
$\theta_*D_{\eta} - D_{\eta} \sim \theta_*D_{0,\eta} - D_{0,\eta}$ on
$X_{\eta}$.  Hence $\theta$ acts on $W(X/S)$ as a translation by
$\frac 1d [\theta_*D_0 - D_0] \in N^1(X/S)/V(X/S)$. 
The map $\sigma': M \to N^1(X/S)/V(X/S)$ given by 
$\sigma'(\theta) = \frac 1d [\theta_*D_0 - D_0]$ is a group homomorphism, 
and its image is a free abelian group.

Let $J_{\eta}$ be the Jacobian of $X_{\eta}$.  
Then $M$ is the group of $\Bbb C(S)$-rational points of $J_{\eta}$.
The natural morphism
$\phi_{D_{0, \eta}}: J_{\eta} \to J_{\eta}$ given by
$\phi_{D_{0, \eta}}(x) = T_x(D_{0,\eta}) - D_{0,\eta}$ is an etale 
morphism of degree 
$d^2$, where $T_x$ denotes the translation by $x$.
For any divisor $D$ on $X$, let $D' = D - (\frac 1d \text{deg }D - 1)D_0$.
Then $\text{deg }D' = d$, and 
$[D'] \in N^1(X/S)/V(X/S)$ is equal to a class of a prime divisor.
If $D$ is a prime divisor of degree $d$, then
the sum of $\phi_{D_{0, \eta}}^{-1}(D_{\eta} - D_{0, \eta})$
in $J_{\eta}$ is an element $x \in M$.
By definition, we have $T_x(D_{0, \eta}) - D_{0, \eta} \sim 
d^2(D_{\eta} - D_{0, \eta})$.
Therefore, the image of $\sigma'$ has rank $\rho(X/S) - v(X/S) -1$.
\hfill $\square$
\enddemo

\proclaim{Theorem 3.6}
Let $f: X \to S$ be a Calabi-Yau fiber space such that $\text{dim }X = 3$
and $\text{dim }S = 2$.
Then there exist only finitely many chambers for the marked minimal models of 
$f$ and finitely many faces of them up to the action of $\text{Bir}(X/S)$,
hence the first parts of Conjecture 1.12 (1) and (2) are true.
\endproclaim

\demo{Proof}
We proceed by induction on $v(X/S) = \text{dim }V(X/S)$.
Assume first that $v(X/S) = 0$.
Since the quotient space $W(X/S)/G(X/S)$ is compact,
there exists only a finite number of faces and chambers up to the action of 
$\text{Bir}(X/S)$ by Theorems 1.9 and 2.6.
We note that there is no fiber space structure in this case.

Next we assume that $v(X/S) > 0$.
Let $I(X/S)$ be the union of all the cones $\alpha_*g_1^*\Cal A(Y'/S)$ for
those marked minimal models $(X', \alpha)$ whose structure morphism
$f': X' \to S$ is factored as $X' @>{g_1}>> Y' @>{g_2}>> T @>h>> S$, where
$h$ is a non-trivial birational morphism.
By the induction hypothesis, there exist 
only finitely many chambers and faces in $I(X/S)$
up to the action of 
$\text{Bir}(X/S)$, because there are only finitely many  
possibilities of $h$ by Lemma 3.3.

Let 
$$
J(X/S) = \{z \in N^1(X/S) ; \text{deg }z = 1 \text{ and } 
z \in \Cal M^e(X/S) \setminus I(X/S)\}.
$$
Then $J(X/S)$ is a closed subset of $N^1(X/S)$ by Theorem 2.6.
Indeed, if $X' \to Y'$ is a birational contraction morphism over $S$ 
such that $Y' \to S$ is not factored by a non-trivial birational morphism
to $S$, and if $Y' \to Z'$ is a further birational contraction over $S$, then
$Z' \to S$ is not factored either.

We shall prove that the natural map $p: J(X/S) \to W(X/S)$ is proper.
Then the theorem follows from Theorems 1.9 and 2.6.
We note that the faces on the boundary of the big cone are already treated in 
Lemma 3.3.
Suppose that there exists a divergent sequence of points $\{z_n\}$ 
in $J(X/S)$ such that $\{p(z_n)\}$ converges in
$W(X/S)$.
We may assume that $f$ factors as $f = h \circ g$ such that $g$ is
equi-dimensional by Lemma 3.3.
Let $C_i$ ($i = 1, \ldots, k$) 
be exceptional divisors of $h$, and $H_i$ divisors on $T$
such that $(H_i \cdot C_i) > 0$ for all $i$ and
$(H_i \cdot C_j) = 0$ for $i \ne j$.
Let $\{y_i\}$ ($i = 1, \ldots, \rho$) for $\rho = \rho(X/S)$ 
be a basis of $N^1(X/S)$ 
such that $y_i = [g^*H_i]$ for $1 \le i \le k$
and $y_i = [D_i]$ for some $g$-exceptional prime divisors $D_i$ for 
$k < i \le v = v(X/S)$.

We write $z_n = \sum_{i=1}^{\rho} a_{ni}y_i$.
By assumption, the sequence $z'_n = \sum_{i= v+1}^{\rho} a_{ni}y_i$
converges to some $z' \in N^1(X/S)$.
If we look at the intersection numbers with the general fibers of the 
$g$-exceptional divisors over their images in $T$, we conclude that
the $a_{ni}$ for $k < i \le v(X/S)$ are bounded, since
$z_n \in \Cal M^e(X/S)$.  
Then by looking at the intersection numbers with the multi-sections of
the morphism $f^{-1}(C_i) \to C_i$ for $1 \le i \le k$, 
we deduce that the $a_{ni}$ for these $i$ are bounded from below.
We may assume that the sequences $\{a_{ni}\}$ 
for $1 \le i \le k'$ diverge to infinity for some $k' \le k$,
and that the sequence $z''_n = \sum_{i= k'+1}^v a_{ni}y_i$
converges to some $z'' \in N^1(X/S)$.
Since $\{z_n\}$ diverges, we have $k' \ge 1$.

Let $f: X @>{g'}>> T' @>{h'}>> S$ be the factorization such that
the exceptional locus of $h'$ consists of the $C_i$ for $1 \le i \le k'$.
We have $[z_n] = [z'_n] + [z''_n]$ in $N^1(X/T')$,
hence $[z' + z''] \in \bar \Cal M(X/T') \cap \Cal B(X/T')$,
since $\text{deg }z_n = 1$.
By using flops over $T'$, we take a minimal model $X'$ of $X$
such that $[z' + z''] \in \Cal A^e(X'/T')$, where we use the same symbol for 
the strict transform of a divisor class by abuse of notation. 
Let $X' @>{g_1}>> Y' @>{g_2}>> T'$ be a factorization such that
$g_1$ is a birational morphism and 
$z' + z'' = g_1^*w$ for a $g_2$-ample class $w \in N^1(Y'/S)$.
We take $N$ a large enough integer such that
$w + N\sum_{i=1}^{k'}g_2^*H_i$ is $h' \circ g_2$-ample.
Since there are only finitely many minimal resolutions of $Y'$,
we may assume that $[z_n] \in \Cal A^e(X'/T')$ for all $n$.
Since $a_{ni} \to \infty$ for $1 \le i \le k'$, we have
$z_n \in \Cal A^e(X'/S)$ for large $n$, a contradiction to the assumption that
$z_n \not\in I(X/S)$.
\hfill $\square$
\enddemo

\proclaim{Corollary 3.7}
Asume that $\rho(X/S) = v(X/S) + 1$.
Then $\bar \Cal A(X/S)$ is a finite rational polyhedral cone,
and there exist only finitely many marked minimal models of $f$.
\endproclaim

\demo{Proof}
Because $\text{Bir}(X/S)$ is finite, there exist only finitely many chambers.
Since the boundary of the big cone is a rational hyperplane, the chambers are 
rational polyhedral.
\hfill $\square$
\enddemo

\definition{Example 3.8}
(1) Let $G$ be a finite group which acts freely in codimension 2
on the product of a K3 or an abelian surface $S$ and an elliptic curve $E$. 
Then $Y = (S \times E)/G$ is $\Bbb Q$-factorial and log terminal.
A minimal model of the projection 
$g: (S \times E)/G \to S/G $ satisfies the condition of Corollary 3.7.

\vskip .5pc

(2) (cf. [R, 6.8]). 
Let $C$ be Kodaira's singular fiber of type $I_2$, i.e., 
$C = C_1 \cup C_2$ is a reduced union of two smooth rational curves which 
intersect transversally at two points $P_1, P_2$.
Let $f: X \to S$ be a versal deformation of $C$ over a germ $(S, Q)$.
$S$ is smooth and 2-dimensional.  $X$ is also smooth.
There are two smooth divisors $\Delta_1, \Delta_2$ on $S$
corresponding to the smoothings of $C$ at the points $P_1, P_2$, respectively.
The fibers of $f$ over $\Delta_1 \cup \Delta_2 \setminus \{Q\}$
are isomorphic to a rational curve with a node.
We have $\rho(X/S) = 2$, and
$\bar \Cal M(X/S) = 
\bar \Cal B(X/S) = \{z \in N^1(X/S); (z \cdot F) \ge 0\}$ for a general
fiber $F$ of $f$.
There are only two curves on $X$ which can be flopped; these are the 
irreducible components $C_1, C_2$ of $C$.
If we flop $C_1$, then we obtain again a versal deformation
$f_1: X_1 \to S$ of a fiber of type $I_2$, and we can flop
the strict transform of $C_2$.
We can flop $X$ infinitely many times in this way, and obtain a decomposition
$$
\Cal M^e(X/S) = \bigcup_{(X', \alpha)} \bar \Cal A(X'/S, \alpha)
= \{z \in N^1(X/S); (z \cdot F) > 0\} \cup \{0\}.
$$
A non-zero class on the boundary of the closed movable cone 
is not realized by a nef divisor on any minimal model, hence
not by an effective divisor.

Since $f$ and $f_1$ are versal deformations of a fiber of type $I_2$, 
they are isomorphic over an isomorphism of $S$.
We claim that they are isomorphic over the identity of $S$.
Indeed, let $\sigma_1$ be a section of $f$ which intersects $C_1$, 
and $\sigma_2$ a birational section of $f$ which is smooth, 
contains $C_1$ and intersects $C_2$ transversally at $P_1$ and $P_2$.
Then the translation of $f^{-1}(S \setminus \{Q\})$ which sends
$\sigma_1 \vert_{S \setminus \{Q\}}$ to 
$\sigma_2 \vert_{S \setminus \{Q\}}$ extends to a birational automorphism of 
$f: X \to S$ which corresponds to the marked minimal model
$f_1: X_1 \to S$, hence $X$ and $X_1$ are isomorphic over $S$. 
Note that the existence of such $\sigma_2$ is a consequence of the
existence of an isomorphism $X \to X_1$ which is not necessarily over $S$.

\vskip .5pc

(3) (cf. [Nk2]). 
More generally, we consider Kodaira's singular fiber $C$ of type $I_{a+b}$
for positive integers $a$ and $b$. 
We have $C = \bigcup_{i=1}^{a+b} C_i$ such that $C_i \cong \Bbb P^1$
with transversal intersections $P_i = C_i \cap C_{i+1}$, where we set
$a + b + 1 = 1$.
Let $\{1, \ldots, a + b \} = I_1 \sqcup I_2$ be a decomposition such that
$\# I_1 = a$ and $\#I_2 = b$.
Let $f: X \to S$ be a deformation of $C$ whose total space $X$ is smooth and 
the base space $S$ is a germ of a smooth surface such that
the discriminant $\Delta = \Delta_1 + \Delta_2$ is a normal crossing divisor 
with 2 irreducible components.
We assume that singular points $P_i$ are smoothed for $i \not\in I_{\epsilon}$ 
and remain singular for $i \in I_{\epsilon}$ along $\Delta_{\epsilon}$ 
for $\epsilon = 1,2$.
Thus the singular fibers of $f$ are of type $I_a$ (resp. $I_b$) over
general points of $\Delta_1$ (resp. $\Delta_2$).

We have $\rho(X/S) = a + b$, and
$\bar \Cal A(X/S)$ is dual to a simplicial cone in $N_1(X/S)$ 
generated by the classes of the curves $C_i$.
A curve $C_i$ can be flopped if $i-1$ and $i$ belong to different
sets among $I_1$ and $I_2$.  
For example, if $i-1 \in I_1$ and $i \in I_2$, then
the flop of $C_i$ yields a new minimal model $f': X' \to S$ which corresponds
to a decomposition given by
$I'_1 = I_1 \cup \{i\} \setminus \{i-1\}$ and
$I'_2 = I_2 \cup \{i-1\} \setminus \{i\}$.
In this way, all the isomorphism classes of the marked minimal models of $X$ 
over $S$ are given by different decompositions of the set
$\{1, \ldots, a + b \}$ into the sets with $a$ and $b$ members up to
the cyclic permutations and reversing the order.
In particular, the number is finite as in (2), but the minimal model is not 
unique if $a, b \ge 2$.

\vskip .5pc

(4) 
Let $V = \Bbb P^2 \times \Bbb P^1 \times \Bbb P^1$, and $X$ a general member
of the linear system $\vert \Cal O_V(3,2,2) \vert$.
Then $K_X \sim 0$, and $\rho(X) = 3$ with generators 
$A = \Cal O_X(1,0,0), B = \Cal O_X(0,1,0)$ and $C = \Cal O_X(0,0,1)$.
The nef cone $\bar \Cal A(X)$ is generated by $A, B$ and $C$.
The faces $OAB$ and $OAC$ correspond to small birational contractions, 
and $OBC$ to an elliptic fibration.
The ray $OA$ corresponds to an elliptic fibration, and
$OB$ and $OC$ to K3 fibrations.

The contraction $\phi: X \to Y$ corresponding to the face $OAB$ is the Stein 
factorization of the projection
$p_{12}: X \to \Bbb P^2 \times \Bbb P^1$.  $\phi$ contracts
$54 = (\Cal O_{\Bbb P^2 \times \Bbb P^1}(3,2)^3)$
$(-1,-1)$-curves.
Let $\phi': X' \to Y$ be its flop.
We denote by $A', B'$ and $C'$ the strict transforms of $A, B$ and $C$,
respectively.
Then $C_1 = 3A' + 2B' - C'$ is nef and $C_1^2 = 0$.  

Let $f(x,y)z_0^2 + g(x,y)z_0z_1 + h(x,y)z_1^2 = 0$ be an equation of $X$,
where $z_0, z_1 \in H^0(\Cal O_V(0,0,1))$ and 
$f, g, h \in H^0(\Cal O_V(3,2,0))$.
Then $fz_1^{-1}$ and $hz_0^{-1}$ give a basis of
$H^0(X, 3A + 2B - C) \cong H^0(X', C_1)$.
So we have an embedding $X' \to V$ with an equation
$f(x,y)z_0^{\prime 2} + g(x,y)z'_0z'_1 + h(x,y)z_1^{\prime 2} = 0$
for $z'_0 = hz_1$ and $z'_1 = fz_0$.
Thus we have an isomorphism $\beta_1: X' \to X$, and 
the nef cone $\bar \Cal A(X')$ is generated by $A', B'$ and $C_1$.
Moreover, the birational automorphism
$\gamma_1 = \beta_1 \circ \alpha_1$ of $X$ for  
$\alpha_1 = \phi^{\prime -1}\circ \phi$ is 
given by $(x,y,z) \mapsto (x,y,fh^{-1}z^{-1})$.
So $\gamma_1^2 = \text{id}$.

Let $\gamma_2 \in \text{Bir}(X)$ be another birational automorphism of $X$ 
corresponding to the face $OAC$.  We set
$C_n = \frac 32 (n^2 + n)A + (n+1)B - nC$ for $n \in \Bbb Z$.
Thus $B = C_0, C = C_{-1}$.
Then the cones $\bar \Cal A(X_n, \alpha_n)$
generated by the $A, C_{n-1}$ and $C_n$ exhaust all the chambers;
$\Cal M^e = \bigcup_n \bar \Cal A(X_n, \alpha_n)$.
Here the marked minimal models $(X_n, \alpha_n)$ correspond to 
the birational automorphisms $\gamma_1 \circ \gamma_2 \circ \gamma_1 \cdots$
of $X$.

All the flops above are relative to the elliptic fibration 
$\pi: X \to \Bbb P^2$ corresponding to the ray $OA$, and
the rays $OC_n$ approach to $OA$ 
when $\vert n \vert \to \infty$.
We note also that $\bar \Cal B(X) = \bar \Cal M(X) = \Cal M^e(X)$, 
because there is no boundary face corresponding to a divisorial contraction.

\vskip .5pc

(5) Let $S = \Bbb P^2$, $\mu: \Sigma_1 \to S$ a blow-up at a point $Q$, and
$g: \Sigma_1 \to \Bbb P^1$ a natural $\Bbb P^1$-bundle.
Let $h: R \to \Bbb P^1$ be a rational elliptic surface with 
infinitely many sections $s_{\lambda}$ which has only singular fibers 
of type $I_1$.
Then we obtain a Calabi-Yau fiber space 
$f = \mu \circ p_1: X = \Sigma_1 \times_{\Bbb P^1} R \to \Sigma_1 \to S$.
A fiber $E = f^{-1}(Q)$ is isomorphic to $R$, any curve in a fiber of $f$
can be deformed to a curve in $E$, and any divisor on $E$ can be extended
to a divisor on $X$ by a morphism $p_2: X \to R$. 
Therefore, we have an isomorphism $p_2^*: N^1(R) \to N^1(X/S)$
such that $p_2^*(\Cal A^e(R)) = \Cal A^e(X/S)$.
$f$ has also infinitely many birational sections 
$D_{\lambda} = p_2^*(s_{\lambda})$, and the translations among them
give infinitely many automorphisms of $f$.
Since the action of $\text{Aut}(R/\Bbb P^1)$ 
on $\Cal A^e(R)$ has a finite rational 
polyhedral fundamental domain ([GM]), so has that of $\text{Aut}(X/S)$
on $\Cal A^e(X/S)$.  Any birational automorphism of $f$ induces that of
$h$.  Since $\text{Bir}(R/\Bbb P^1) = \text{Aut}(R/\Bbb P^1)$, 
we have $\text{Bir}(X/S) = \text{Aut}(X/S)$.

The intersections $C_{\lambda} = D_{\lambda} \cap E$ are $(-1,-1)$-curves,
and we can flop some of them to obtain new minimal models of $f$.
Let $f': X' \to S$ with $\alpha: X' -\to X$ 
be a marked minimal model obtained by flopping 9 disjoint
curves among the $C_{\lambda}$.
Then the strict transform $E'$ of $E$ is isomorphic to $\Bbb P^2$.
Since any curve in a fiber of $f'$ can again be deformed to a curve in 
$f^{\prime -1}(Q)$, the nef cone $\bar \Cal A(X'/S)$
itself is a finite rational polyhedral cone.
Since there are infinitely many choices of the 9 disjoint curves, we obtain
infinitely many different chambers.
But the number of sets of disjoint 9 sections of $h$ 
up to $\text{Aut}(R/\Bbb P^1)$ is finite, 
because the number of sections which are disjoint from a given section 
is finite.
Since a biregular automorphism of $f$ induces a birational automorphism of 
$f'$, the number of chambers up to $\text{Bir}(X'/S)$ is finite. 
\enddefinition

\head 4.  K3, abelian, Enriques or hyperelliptic fiber spaces
\endhead

We shall consider a Calabi-Yau fiber space $f: X \to S$ such that
$\text{dim }X = 3$ and $\text{dim }S = 1$ in this section,
and prove the first parts of Conjecture 1.12 (1) and (2) in this case.
Let $X_{\eta}$ denote the generic fiber, 
and $X_{\bar \eta}$ the geometric generic fiber for the algebraic closure
$\overline{\Bbb C(S)}$ of $\Bbb C(S)$.
A prime divisor $D$ on $X$ is called {\it $f$-vertical} if $f(D) \ne S$. 
Moreover, if $[D] \ne 0$ in $N^1(X/S)$, then $D$ is $f$-exceptional.
There are $f$-exceptional divisors which are not $f$-vertical.
Let $V(X/S)$ be the subspace of $N^1(X/S)$ generated by the classes of
$f$-vertical divisors, and $v(X/S) = \text{dim }V(X/S)$.
$V(X/S)$ is stable under the action of $\text{Bir}(X/S)$, and 
there exists a subgroup of finite index of $\text{Bir}(X/S)$ which acts
trivially on $V$.
By Theorem 2.1 and Remark 2.2, 
there exists a rational polyhedral fundamental domain
$\Pi(X_{\eta})$ for the action of $\text{Aut}(X_{\eta}) \cong \text{Bir}(X/S)$
on $\Cal A^e(X_{\eta})$.

\proclaim{Lemma 4.1}
(1) The natural restriction homomorphism
$r: N^1(X/S) \to N^1(X_{\eta})$
defined by $r([D]) = [D_{\eta}]$ is surjective.

(2) $r(\Cal M^e(X/S)) = \Cal A^e(X_{\eta})$, 
$r^{-1}(\Cal A^e(X_{\eta})) \cap \bar \Cal M(X/S) = \Cal M^e(X/S)$, and
\break
$r^{-1}(\Cal B(X_{\eta})) = \Cal B(X/S)$.

(3) Let $r_1: N^1(X/S) \to N^1(X/S)/V(X/S)$ be the projection.
Then its restriction to $\bar \Cal M(X/S)$ is proper. 

(4) The natural homomorphism 
$s: N^1(X_{\eta}) \to  N^1(X_{\bar \eta})$ defined by
$s([D_{\eta}]) = [D_{\bar \eta}]$ is injective,
and its image coincides with the invariant part
$N^1(X_{\bar \eta})^G$ for $G = \text{Gal}(\overline{\Bbb C(S)}/\Bbb C(S))$.
\endproclaim

\demo{Proof}
(1) Since any prime divisor on $X_{\eta}$ can be extended to $X$
as a $\Bbb Q$-Cartier divisor, $r$ is surjective.

(2) Since a movable divisor on a surface is nef, 
we have $r(\Cal M^e(X/S)) \subset \Cal A^e(X_{\eta})$. 
If $D_{\eta}$ is an effective and nef divisor on $X_{\eta}$, then
some positive multiple $mD_{\eta}$ is free, hence its closure in $X$
is movable, and
$r(\Cal M^e(X/S)) \supset \Cal A^e(X_{\eta})$. 

If $D_{\eta}$ is effective, then $D$ is $f$-effective, hence the second 
assertion.  The third assertion is obvious.

(3) Let $F_i = \sum_{j=1}^{c_i} m_{ij}D_{ij}$ for $i = 1, \ldots, 
d$ be the reducible fibers of $f$, and let $H$ be an $f$-ample divisor. 
Then $v(X/S) = \text{dim }V(X/S) = \sum_i (c_i - 1)$, and
the $D_{ij}$ for $1 \le i \le d$ and $1 \le j < c_i$ form a basis of $V(X/S)$.
If a divisor $D$ is movable, then $(H \cdot D \cdot D_{ij}) \ge 0$ for any
$i, j$.  Since the matrix $[(H \cdot D_{ij} \cdot D_{k\ell})]
_{1 \le i, k \le d, 1 \le j < c_i, 1 \le \ell < c_k}$
of size $v(X/S)$ is negative definite, we have our assertion.

(4) is clear.
\hfill $\square$
\enddemo

We fix an ample class $\xi$ on $X_{\eta}$.
We consider an affine subspace  
$N^1_1(X_{\eta}) = \{w \in N^1(X_{\eta}) ; (w \cdot \xi) = 1 \}$ of
$N^1(X_{\eta})$.  Similarly, we define
$N^1_1(X/S) = r^{-1}(N^1_1(X_{\eta}))$, 
$\Pi_1(X_{\eta}) = \Pi(X_{\eta}) \cap N^1_1(X_{\eta})$, and so on.
Let
$$
W(X/S) = \{[z] \in N^1_1(X/S)/V(X/S) ;  
r(z) \in \Pi_1(X_{\eta}) \cap \Cal B(X_{\eta}) \}.
$$

\proclaim{Lemma 4.2}
If $X_{\eta}$ is a K3 or Enriques surface, then
the kernel of $r$ coincides with the subspace $V(X/S)$, hence
$\rho(X/S) = \rho(X_{\eta}) + v(X/S)$
and $W(X/S) \cong \Pi_1(X_{\eta}) \cap \Cal B(X_{\eta})$.
\endproclaim

\demo{Proof}
Let $D$ be a divisor on $X$ such that $D_{\eta} \equiv 0$.
Then $2D_{\eta} \sim 0$, and $2D$ is vertical.
\hfill $\square$
\enddemo

\proclaim{Lemma 4.3}
Assume that $X_{\bar \eta}$ is an abelian surface  
or a hyperelliptic surface.
Then the image of the representation
$\sigma: \text{Bir}(X/S) \to GL(N^1(X/S), \Bbb Z)$
contains a subgroup $G(X/S)$ which satisfies the following conditions:

(1) $G(X/S)$ acts trivially on $N^1(X_{\eta})$ and $V(X/S)$.

(2) $G(X/S)$ is isomorphic to a free abelian group of rank
$\rho(X/S) - \rho(X_{\eta}) - v(X/S)$.

(3) $G(X/S)$ acts on the fibers of the projection 
$W(X/S) \to \Pi_1(X_{\eta})$
properly discontinuously as a group of translations, and
the quotient space $W(X/S)/G(X/S)$ is a real 
torus bundle over $\Pi_1(X_{\eta})$.
\endproclaim

In the case of Lemma 4.2, we set $G(X/S) = \{\text{id}\}$.

\demo{Proof}
Instead of (3), noting that
$\Cal A^e(X_{\eta}) \cap \Cal B(X_{\eta}) = \Cal A(X_{\eta})$,
we shall prove that $G(X/S)$ acts on the fibers of the projection 
$$
\{[z] \in N^1(X/S)/V(X/S) ; r(z) \in \Cal A(X_{\eta}) \} \to \Cal A(X_{\eta})
$$
properly discontinuously as a group of translations, and
the quotient space is a real torus bundles over $\Cal A(X_{\eta})$.

First, assume that $X_{\bar \eta}$ is an abelian surface. 
Let $D_0$ be a divisor such that $D_{0, \eta}$ is ample.
Let $M$ be the subgroup of $\text{Bir}(X/S)$ consisting of elements
which act on $X_{\eta}$ as translations and on $V(X/S)$ trivially.
If we put $G(X/S) = \sigma(M)$, then we have (1).
For $\theta \in M$ and for a divisor $D$ such that 
$[D_{\eta}] = [D_{0, \eta}]$ in $N^1(X_{\eta})$, 
we have $\theta_*D_{\eta} - D_{\eta} \sim \theta_*D_{0,\eta} - D_{0,\eta}$ on
$X_{\eta}$.  
The map $\sigma': M \to N^1(X/S)/V(X/S)$ given by 
$\sigma'(\theta) = [\theta_*D_0 - D_0]$ is a group homomorphism, and its image
is a free abelian group.

$X_{\eta}$ is a torsor of an abelian surface $A_{\eta}$,
and $M$ is the group of $\Bbb C(S)$-rational points of $A_{\eta}$.
The natural morphism
$\phi_{D_{0, \eta}}: A_{\eta} \to A^*_{\eta}$ given by
$\phi_{D_{0, \eta}}(x) = T_x(D_{0,\eta}) - D_{0,\eta}$ is an etale 
morphism, where $A^*_{\eta}$ is the dual abelian surface and
$T_x$ denotes the translation by $x$.
If $D$ is a divisor such that $[D_{\eta}] = [D_{0, \eta}]$ in
$N^1(X_{\eta})$, 
the sum of $\phi_{D_{0, \eta}}^{-1}(D_{\eta} - D_{0, \eta})$
in $A_{\eta}$ gives an element $\theta = T_x \in M$ such that
$T_x(D_{0, \eta}) - D_{0, \eta} \sim d(D_{\eta} - D_{0, \eta})$ 
for $d = \text{deg }\phi_{D_{0, \eta}}$.
Therefore, the image of $\sigma'$ has rank 
$\rho(X/S) - \rho(X_{\eta}) - v(X/S)$, and have (2) and (3).

Next, assume that $X_{\bar \eta}$ is a hyperelliptic surface.
Then there exists a cyclic Galois coverring 
$\tau: \tilde X_{\eta} \to X_{\eta}$ such that
$\tilde X_{\eta}$ is a torsor over an abelian surface $\tilde A_{\eta}$.
Let $\zeta$ be a genrator of $\text{Gal}(\tilde X_{\eta}/X_{\eta})$.
Let $D_0$ and $D$ be divisors on $X$ such thath
$[D_{0, \eta}] = [D_{\eta}] \in \Cal A(X_{\eta})$.
Then we have a translation $\theta = T_x$ of $\tilde X_{\eta}$
such that $\theta_*\tau^*D_{0, \eta} - \tau^*D_{0, \eta}
\sim d(\tau^*D_{\eta} - \tau^*D_{0, \eta})$ for some $d$
by the previous argument.
Since $x \in \tilde A_{\eta}$ is $\zeta$-invariant,
we have $\theta \circ \zeta = \zeta \circ \theta$, 
and $\theta$ descends to an element of 
$\text{Aut}(X_{\eta}) \cong \text{Bir}(X/S)$.
Now the proof is reduced to the previous case.
\hfill $\square$
\enddemo

\proclaim{Theorem 4.4}
Let $f: X \to S$ be a Calabi-Yau fiber space such that $\text{dim }X = 3$
and $\text{dim }S = 1$.
Then there exist finitely many chambers and faces up to the action of 
$\text{Bir}(X/S)$ for the marked minimal models of $f$.
\endproclaim

\demo{Proof}
First, we consider the chambers for the marked minimal models of $f$
which have faces corresponding to fiber space structures.  
Two fiber space structures 
$X_1 @>{g_1}>> T_1 @>{h_1}>> S$ and $X_2 @>{g_2}>> T_2 @>{h_2}>> S$ are
said to be birationally equivalent if 
there exist birational maps $\alpha: X_1 -\to X_2$ and
$\beta: T_1 -\to T_2$ over $S$ such that $\beta \circ g_1 = g_2 \circ \alpha$.
By Theorem 2.1 and Remark 2.2, 
the generic fiber $X_{\eta}$ has only finitely many 
elliptic surface structures up to biregular automorphisms.
Since any biregular automorphism of $X_{\eta}$ extends to a birational 
automorphism of $f: X \to S$, there are only finitely many elliptic fiber 
space structures up to birational equivalence.
By Lemma 3.3, we conclude that there exist only finitely many faces 
corresponding to fiber space structures up to $\text{Bir}(X/S)$.
Then by Theorem 3.6, we deduce our finiteness assertion for the chambers 
and their faces which in turn have faces corresponding to 
fiber space structures, because a birational automorphism
of $X$ over $T$ is also over $S$. 

Let $I(X/S)$ be the union of all the cones $\alpha_*g_1^*\Cal A(Y'/S)$ for
those marked minimal models $(X', \alpha)$ whose structure morphism
$f': X' \to S$ is factored as $X' @>{g_1}>> Y' @>{g_2}>> T @>h>> S$, where
$\text{dim }T = 2$.
The above argument showed that there exist 
only finitely many chambers and faces in $I(X/S)$
up to the action of $\text{Bir}(X/S)$.

Next we consider the interior of the big cone $\Cal B(X/S)$.
Let $\tilde \Pi_1(X/S) = 
\mathbreak
r^{-1}(\Pi_1(X_{\eta})) \cap \Cal M^e(X/S)$.
In the rest of the proof, we shall prove that
$$
J(X/S) = \tilde \Pi_1(X/S) \setminus I(X/S)
$$
is a closed subset contained in $\Cal B_1(X/S)$.
Then the projection $J(X/S) \to 
\mathbreak
N_1^1(X/S)/V(X/S)$ is proper, and
the theorem follows from Lemmas 4.2, 4.3 and Theorems 1.9, 2.6.
Since $J(X/S) \cap \Cal B(X/S)$ is a closed subset of $\Cal B(X/S)$
as in the proof of Theorem 3.6,
it is sufficient to prove that there exists no
element $z$ in the closure of $J(X/S)$ such that $w^2 = 0$ for $w = r(z)$.
Assuming the contrary, we shall derive a contradiction.

We claim that $G(X/S) = \{\text{id}\}$ in this case.
Indeed, under our assumption that there exists a rational point $w$, 
there exist points $w_i \in \Cal A^e(X_{\eta})$ which generate the
vector space $N^1(X_{\eta})$ and such that $w_i^2 = 0$.  
Let $X_i @>g_i>> T_i @>h_i>> S$ be factorizations of some marked minimal models
of $f$ such that all the fibers of the $g_i$ are purely 1-dimensional 
and $w_i = (g_i^*v_i)_{\eta}$ for some $h_i$-ample classes $v_i$ on the $T_i$.
Since $G(X/S)$ acts trivially on $N^1(X_{\eta})$, it induces automorphisms
on the $T_i$, and there exists a subgroup of finite index
of $G(X/S)$ which acts trivially on the $v_i \in N^1(T_i/S)$.
Since $G(X/S)$ acts on the fibers of $r$ as groups of translations,
we conclude that the action of $G(X/S)$ on $N^1(X/S)$ is trivial.

We may assume that there exists a factorization $X @>g>> T @>h>> S$ of $f$ 
such that $g$ is equi-dimensional.
Let $\{x_i\}$ ($i = 1, \ldots, \rho$) for 
$\rho = \rho(X/S) = \rho(X_{\eta}) + v(X/S)$ be a basis of $N^1(X/S)$ 
such that the $r(x_i)$ for $1 \le i \le \rho_{\eta} = \rho(X_{\eta})$ 
form a basis of
$N^1(X_{\eta})$ with $w = r(x_{\rho_{\eta}})$, and that
the $x_i$ for $\rho_{\eta} < i \le \rho(X/S)$ form a basis of $V(X/S)$.
We also assume that there exists a basis $\{y_j\}$ ($j =1, \ldots, \rho(T/S)$) 
of $N^1(T/S)$ such that $g^*y_1 = x_{\rho_{\eta}}$, 
the $y_j$ for $j > 1$ are classes of $h$-exceptional prime divisors, 
$g^*y_j = x_{\rho_{\eta}+j-1}$ for $j > 1$, and that
the $x_i$ for $i \ge \rho_{\eta} + \rho(T/S)$ 
are classes of $g$-exceptional prime divisors.

Let us take a sequence $\{z_n\}$ in $J(X/S)$ which converges to $z$, and
write $z_n = \sum_{i=1}^{\rho} a_{ni}x_i$.
Then we can take a sequence of positive numbers $\{c_n\}$ which diverges to 
infinity such that the sequence 
$z'_n = \sum_{i= 1}^{\rho_{\eta}-1} a_{ni}c_nx_i$
converges to some non-zero $z' \in N^1(X/S)$ if we replace the sequence
$\{z_n\}$ by its suitable subsequence.
If we look at the intersection numbers with the general fibers of the 
$g$-exceptional divisors over their images in $T$, we conclude that
the $a_{ni}c_n$ for $i \ge \rho_{\eta} + \rho(T/S)$ are bounded, since
$z_n \in \Cal M^e(X/S)$. 
Therefore, we may assume that 
$z''_n = \sum_{i= \rho_{\eta} + \rho(T/S)}^{\rho} a_{ni}c_nx_i$
converges to some $z'' \in N^1(X/S)$

We take generators $v_k$ ($k = 1, \ldots, e$) 
of the rational polyhedral cone $\Cal A^e(T/S)$.
Then we can write 
$\sum_{j= 1}^{\rho} a_{n,\rho_{\eta}+j-1}c_ny_j
= \sum_{k=1}^e b_{nk}v_k$ for some numbers $b_{nk}$.
We note that these expressions are not unique,
but we may assume that the $b_{nk}$ are bounded from below.
Then we may assume that the sequences $\{b_{nk}\}$ 
for $1 \le k \le e'$ diverge to infinity for some integer $e' \le e$
and the sequence $y'_n = \sum_{k=e'+1}^e b_{nk}v_k$ converges to some
$y' \in N^1(T/S)$.
Since the sequence $\{c_n\}$ diverges, we have $e' \ge 1$ and
$\sum_{k=1}^{e'} v_k$ is $h$-big.

Let $T @>{h_1}>> T' @>{h_2}>> S$ be the factorization of $h$ such that
$\sum_{k=1}^{e'} v_k = h_1^*d$ for some $h_2$-ample class $d \in N^1(T'/S)$.
We have $[c_nz_n] = [z'_n] + [z''_n] + [g^*y'_n]$ in $N^1(X/T')$,
hence $[z' + z'' + z'''] \in \Cal M^e(X/T')$, where $z''' = g^*y'$.
Since $z_n \in \Cal M^e(X/S)$ and $r(z_n) \in \Pi(X_{\eta})$, we have 
$[z' + z'' + z'''] \in \bar \Cal M(X/T') \cap \Cal B(X/T')$.
By using flops over $T'$, we take a minimal model $X'$ of $X$
such that $[z' + z'' + z'''] \in \Cal A^e(X'/T')$.
Let $X' @>{g_1}>> Y' @>{g_2}>> T'$ be a factorization 
such that $g_1$ is a birational morphism and
$z' + z'' + z''' = g_1^*\ell$ for a $g_2$-ample class $\ell \in N^1(Y'/S)$.
We take $N$ a large enough integer such that
$\ell + Ng_2^*d$ is $h_2 \circ g_2$-ample.
Since there are only finitely many minimal resolutions of $Y'$,
we may assume that $[z_n] \in \Cal A^e(X'/T')$ for all $n$.
Since $b_{nk} \to \infty$ for $1 \le k \le e'$, we have
$z_n \in \Cal A^e(X'/S)$ for large $n$, a contradiction to the assumption that
$z_n \not\in I(X/S)$.
\hfill $\square$
\enddemo

As a corollary to Theorems 3.6 and 4.4, we obtain

\proclaim{Theorem 4.5}
Let $X$ be an algebraic variety of dimension 3 whose Kodaira dimension 
$\kappa (X)$ is positive.
Then there exist only finitely many minimal models of $X$
up to isomorphisms.
\hfill $\square$
\endproclaim

\definition{Remark 4.6}
The argument of the above proof can be generalized 
to investigate the local structure of chambers near a face
corresponding to a fiber space structure in general.
\enddefinition

\enddocument